# MATHEMATICAL ASPECTS OF THE PERIODIC LAW

**Guillermo Restrepo*§ and Leonardo A. Pachón**¶**

*Laboratorio de Química Teórica, Universidad de Pamplona,
Pamplona, Norte de Santander, Colombia.
** Escuela de Física, Facultad de Ciencias,
Universidad Industrial de Santander, Santander, Colombia.

**Abstract**:
We review different studies of the Periodic Law and the set of chemical elements from a mathematical point of view. This discussion covers the first attempts made in the 19th century up to the present day. Mathematics employed to study the periodic system includes number theory, information theory, order theory, set theory and topology. Each theory used shows that it is possible to provide the Periodic Law with a mathematical structure. We also show that it is possible to study the chemical elements taking advantage of their phenomenological properties, and that it is not always necessary to reduce the concept of chemical elements to the quantum atomic concept to be able to find interpretations for the Periodic Law. Finally, a connection is noted between the lengths of the periods of the Periodic Law and the philosophical Pythagorean doctrine.

**Key words**
Chemical Elements, Periodic Law, Number Theory, Information Theory, Order Theory, Set Theory, Topology, Mathematical Chemistry

## 1. INTRODUCTION

One of the most important constructs of chemical knowledge is the Periodic Table. This graphical representation of the underlying law of the chemical elements, the Periodic Law, has been studied from several points of view, the majority of them related to physics (Scerri et al., 1998). These studies are focused on the atom as a refinement of the concept of element to try to explain the similarities and trends among the chemical elements, which is what at the end, the Periodic Law shows. To do this, scientists typically defer to quantum chemistry with the aim of explaining the features of the Periodic Law and try to *predict* new results. In general, this kind of work in science is a relevant part of the scientific activity, but another important part concerns the organisation or *description* of known results. In other words, science has two main fields of work, the descriptive and the inferential ones. Despite the inferential character of contemporary science, the Periodic Law was the product of a descriptive analysis of the chemical elements. This was the work done by Mendeleev in 1869 finding patterns in the properties of the chemical elements. In this way, the Russian scientist, trying to complete his description of the chemical elements, postulated the Periodic Law and a graphical representation of it, the Periodic Table. The work of Mendeleev is a special case of descriptive science (Villaveces, 2000). However, he opened the way to cross the hypothetical border between inferential and descriptive science when he was able to *predict* properties of unknown substances[1] (Greenwood and Earnshaw, 2002). All his predictions were based on the description that he made regarding the chemical elements. What Mendeleev really found was a general Periodic Law including a way to make predictions starting from the description. In other words, the Mendeleevian methodology can be synthesised as 'from description to inference.' But the descriptive research of the chemical elements is little subsequently considered in a formal manner. In spite of presently having more than 700 periodic tables (Mazurs, 1974), these are not new descriptions of the chemical elements or their behaviour. They are just different representations of the same phenomena, different shadows of the same object - the Periodic Law. Thence all the periodic tables that we have nowadays are *meta*-descriptions of the fundamental description made by Mendeleev. But, do we have descriptive studies that do not take Mendeleev's description as their basis? The answer is yes! and we discuss some of the research done in this direction in this paper.

§E-mail: grestrepo@unipamplona.edu.co

¶ E-mail: leaupaco@tux.uis.edu.co

On the other hand, as we mentioned above, the majority of the contemporary studies of the chemical elements have been taken from physics. But an important fact that was shown by Mendeleev in his seminal work (Mendeleev, 1869) was the periodic nature of the properties of chemical elements, all of them considered as a whole. In other words, Mendeleev showed that the properties of the chemical elements have a mathematical structure[2]. If we consider an element as a collection of its properties, then we can say that the Russian scientist showed that there is a mathematical structure in the set of chemical elements. We describe here different mathematical approaches to the Periodic Law, some of them related to the concept of order, of information theory, of similarities, and of topology. Our discussion includes all these attempts which are evidences of different mathematical structures of the chemical elements.

## 2. FROM TWO-DIMENSIONAL PLOTS TO QUANTUM MECHANISCS

Before mentioning the studies done in the search for the mathematical structure of the Periodic Law, it is important to say something regarding the work of Mendeleev and of Meyer. These two scientists found (Mendeleev, 1869; Meyer, 1870) that when one plots a physico-chemical property against the atomic weight (or in modern terms, atomic number Z (Moseley, 1913)), we find an oscillating plot. Mendeleev called such oscillations 'periods' due to the similarity between such graphics to periodic ones of the trigonometric functions, such as sine and cosine. Also, such plots are called periodic by chemists due to the fact that each oscillation covers a number of elements if they are organised according to their atomic number. Thus, we find for the first seven oscillations (or periods) the following cardinalities[3] 2, 8, 8, 18, 18, 32 and 32. It is evident that they are not equal and the periods are not strictly periodic[4] (Babaev and Hefferlin, 1996) but there is a symmetry in their distribution that makes us suppose an underlying mathematical structure for the properties of chemical elements (Villaveces, 2000).

This mathematical sense of the Periodic Law is not a novel discovery. Mendeleev, in his Faraday lecture in 1889, (Mendeleev, 1889) offered a review of the mathematical approaches developed up to his days trying to find the mathematical nature of the chemical periodicity. In that lecture, the Russian scientist mentioned the work of Mills (1886) (Mendeleev, 1889), who considered that all atomic weights could be expressed by the following function: $15(n - 0.9375^t)$, where $n$ and $t$ are *integer numbers*. In that work, oxygen had $n$=2 and $t$=1, which produced an atomic weight of 15.94. Other examples were chloride, bromide and iodide for which $n$ took values of 3, 6 and 9, respectively, while $t$=14, 18 and 20. It is important to say that Mills tried to express the atomic weight into mathematical terms but not the whole Periodic Law. However, in the case of having a mathematical expression for the atomic weight, it would be reasonable to think about the mathematisation of the Periodic Law, due to the capital importance of the atomic weight to the Periodic Law[5]. Mendeleev also discussed the work of Tchitchérin (1888) (Mendeleev, 1889), who actually studied the alkali metals. This scientist found simple relationships among atomic volumes at 25°C, which could be expressed according to $A(2 - 0.00535 A \times n)$, where $A$ is the atomic weight and $n$ is an *integer number*. For example $n=8$ for lithium and sodium, 4 for potassium, 3 for rubidium, and 2 for cesium. These two researches show the early importance of the integer numbers in mathematical studies of the chemical elements. However, the most known and far more fundamental work regarding integers in the Periodic Law was developed by Moseley (Moseley, 1913) years later.

On the other hand, Lewis (Lewis, 1916), studied the periodic system through the chemical properties of the compounds (Cruz et al., 1991) and he proposed an atomic model based on the concept of spatial arrangements of electrons (shells), the first shell with two electrons and thereafter with eight to be placed in the vertices of a cube[6] (Lewis, 1916). Subsequently, Bohr (Bohr, 1922), making use of quantum theory, related the regularities in the behaviour of chemical elements to a positive *integer number n*, today known as the principal quantum number. Afterwards, Madelung[7] (Madelung, 1936) suggested a way (Madelung, Aufbau or Bohr Rule) to construct the periodic table taking advantage of $n$ and the angular momentum (azimutal, quantum number $l$). An 'explanation' of the Periodic Law on the basis of this rule has become a standard cornerstone, so that many scientists have studied the nature of the Madelung Rule making use of different quantum theoretic methodologies[8]. However, these attempts do not focus on the overall mathematical nature of the Periodic Law but rather on the quantum-mechanically mediated electronic structure of the various atoms.

### 3. INFORMATION THEORY AND ORDER

In the 1970s Bonchev et al. applied Shannon's information theory to characterise the structure of atoms and atomic nuclei (Bonchev, 2006). According to Bonchev, "the culmination of the work during that period was the generalization that the electron distribution in atoms obeys a certain principle of maximum information content" (Bonchev, 2006). In spite of the success in the prediction of some properties[9], these studies do not consider the chemical elements from their experimental properties but are an example of a quantum-theoretic approach based on atomic electronic structure. This situation has occurred since 1913 once Moseley linked the concept of element with atomic number as determined by atomic nuclei. On the other hand, Klein (Klein, 1995, 2000; Klein and Babić, 1997) carried out a study of the mathematical aspects of chemical elements defining the set of chemical elements as a multi-partially ordered set (poset). Klein took advantage of the order relations extracted by Mendeleev which were shown in a conventional Periodic Table (Fernelius and Powell, 1982). According to the author (Klein, 1995) the arrangement of the Periodic Table into columns and rows is an evidence of partial orderings[10] of the chemical elements. Klein says (Klein, 1995) "the ionization potentials of elements arranged in a suitable typical periodic chart generally decrease in proceeding down columns and in proceeding right-to-left across rows." This suggests that the Periodic Table may be considered as a multipartially ordered set[11]. Earlier Kreinovich et al. (Scerri et al., 1998) developed a mathematical study of the Periodic Law based on the concept of total order applied to the Madelung rule[12]. Kreinovich et al. suggest that the order between two shells depends on how these shells are related to each other and does not depend on how many shells are hidden inside[13]. This means, according to the authors, that "the two shells $(n,l)$ and $(n',l')$ should depend only on the difference" $n$-$n'$ and $l$-$l'$ (considering $(n,l) < (n',l')$) but not on the absolute values of $n$ and $l$. This property of the Madelung rule is what Kreinovich et al. call a 'local order[14].' Then, the authors show that this local order is a property of Madelung's order ($1s<2s<2p<3s<3p<4s<3d<4p<5s<4d<5p<6s<4f<5d<...$) and also the hydrogenic order or that of ions ($1s<2s<2p<3s<3p<3d<4s<4p<4d<4f<5s<5p<5d<...$)[15]. Kreinovich et al. finally show that it is possible to have "a 1-parametric family of orders that describes the transition from Bohr's order [Madelung one] (for neutral atoms) to the hydrogenic order (that describes highly ionized atoms)" (Scerri et al., 1998). This parameter is the slope $k$ of a ray $r$ that can be defined as $r = \{(x, k \cdot x) \mid x > 0\}$. We show the (infinite) number of $r$'s (shaded region) that can be found in between the Madelung's order ($k$=-1) and the hydrogenic order ($k \to -\infty$) in figure 1.

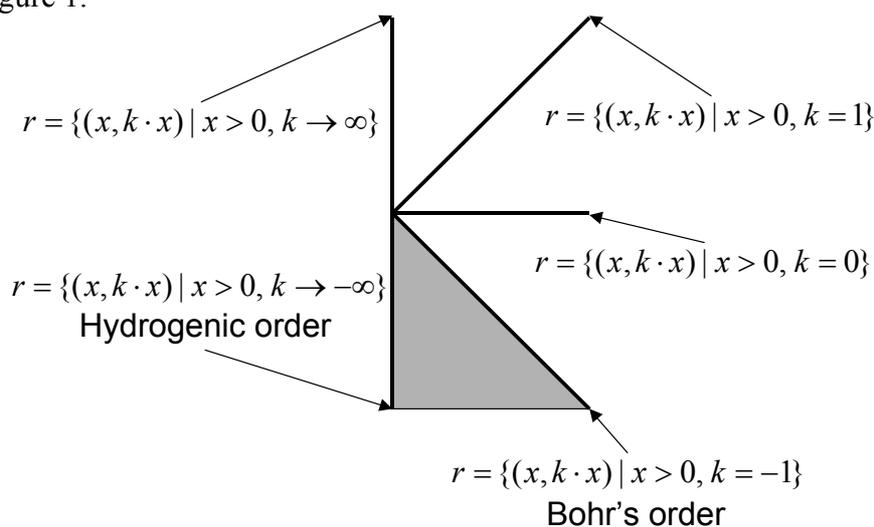

**Figure 1.** Kreinovich's plot of Madelung's order and the hydrogenic order.

Another study related to the order of the chemical elements, this one taking advantage of their properties and the properties of their compounds, was made by Pettifor (Pettifor, 1984, 1986, 1988a, 1988b, 1996)[16]. These studies were made to try to find a map of crystal structures. Pettifor introduced a phenomenological parameter to characterise each element. By means of this he obtained a structural separation of the compounds into two-dimensional plots (Sutton, 1996). These plots offer a separation of the structures only if the elements are ordered along each axis of the plot in a different order than they appear in the traditional Periodic Table (Fernelius and Powell, 1982). This means that the order criterion is not the atomic number but the phenomenological Pettifor order[17].

## 4. SIMILARITY AND CHEMOTOPOLOGICAL STUDIES

Another point of view of the Periodic Law and the Mendeleev methodology is from the set theory. What Mendeleev did was to build up a set $Q$ of chemical elements $q$. Thus, we can define $Q$ as $Q$={$q \mid q$ is a chemical element}. But, how can we define $q$? According to Mendeleev $q$ is determined by its properties[18]. Then we can say that $q$={$x_i \mid x_i$ is the value of $i$-th property for $q$}. Once defining $Q$ and $q$, Mendeleev looked for similarities in the atomic weights and realised the periodicity of some properties as the atomicity (valence). This was possible when he compared the periods of $Q$ to each other. But what Mendeleev did when he classified the elements taking into consideration their atomic weights and their properties, was to partition $Q$. We show a schematic partition of $Q$ from the point of view of the atomic weight (Figure 2a) and then, according to their properties (Figure 2b).

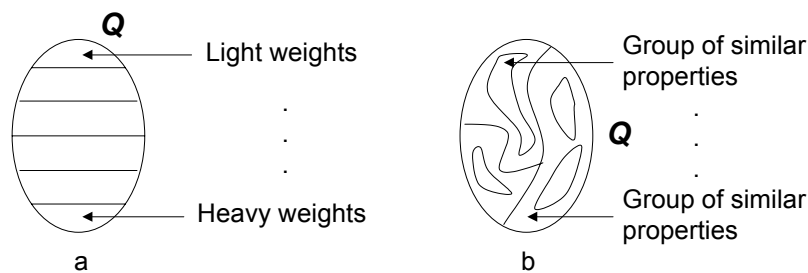

**Figure 2.** Partition of *Q* according to the a) atomic weight and b) the properties of the elements in *Q*.

However in Mendeleev's time set theory was only developing, achieving a consolidated axiomatisation in Frege's book of 1893 (Frege, 1893) and 1903 (Frege, 1903). Recently, some authors have retaken Mendeleev's methodology of studying the elements from the similarity among their properties (Robert and Carbó-Dorca, 1998, 2000; Zhou et al., 2000; Sneath, 2000, Restrepo et al., 2004a. 2004b, 2006a, 2006b). The first attempt in this respect was made by Carbó-Dorca et al. (Robert and Carbó-Dorca, 1998, 2000). These authors used the concept of atomic quantum similarity as an application of the concept of molecular quantum similarity, developed earlier (Carbó et al., 1980, Carbó and Domingo, 1987). These authors studied 20 atoms (Robert and Carbó-Dorca, 1998, 2000) {*Be,C,N,O,Ne,Mg,S,Ar,Ca,Ti,Cr,Ni,Ge,Kr,Zr,Mo,Ru,Cd,Sn,Te*} considering that the information of the system was contained in the electron density function of each one of the 20 atoms. They compared through a quantum similarity measure the similarities among the 20 electron density functions. One of the most important results was that the self-similarity measure calculated for each atom gives a measure of spatial occupation of matter (Robert and Carbó-Dorca, 1998, 2000). However, this study does not consider several atoms of the chemical elements, so that it did not establish trends to find relationships with patterns of the chemical elements.

Other attempts sought similarities among the chemical elements using cluster analysis as a mathematical tool to find the relationships (Zhou et al., 2000; Sneath, 2000, Restrepo et al., 2004a, 2004b, 2005a, 2006a, 2006b). In general, cluster analysis looks for similarities among a set of objects (Aldendefer and Blashfield, 1984; Gordon, 1981; Everitt, 1995). The notion of similarity is normally introduced through a distance or metric function. In this way it is said that *a* and *b* are similar if they are close to each other in the mathematical space where they are defined[19]. Cluster analysis leads to a graph called dendrogram that may be interpreted as a map of similarities (Restrepo et al., 2006c). For instance, in figure 3 a hypothetical dendrogram of 4 elements appears. This dendrogram shows that *a* and *b* are the most similar elements in the set *X*={*a,b,c,d*}. The element *c* is similar to the couple {*a,b*} and finally *d* is similar to the cluster {*a,b,c*}.

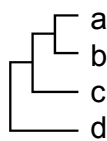

**Figure 3.** A dendrogram of 4 elements.

Keeping in mind the Mendeleev methodology of looking for similarities based on the properties of the elements, some authors combined this idea with cluster analysis to study the chemical elements (Zhou et al., 2000; Sneath, 2000, Restrepo et al., 2004a, 2004b, 2006a, 2006b). Zhou et al. (Zhou et al., 2000), who applied fuzzy cluster analysis[20], developed the first of these approaches. These authors studied a set of 50 elements (Z=1-50) defining each element by 7 physical properties. Some of the clusters found were: {*Co,Ni,Fe,Rh,Ru*}, {*Mo,Tc*}, {*Sc,Y,Ti*}, {*Ga,In,Sn*}, {*N,O,F,H*}, {*Cl,Br*}, {*Zn,Cd*}, {*Ar,Kr,Ne,He*}, {*Mg,Ca*}, {*Li,Na,K,Rb*}, {*Cu,Ag*}, and {*B,C*}. It is important to note the presence of well-known sets of elements such as the noble gases or alkali metals. Afterwards, Sneath[21] applied cluster analysis to study a set of 69 elements (Z=1-83, omitting Z=58-71[22]) using 54 physical and chemical properties (Sneath, 2000). Some of the groups found were: {*He,Ne,Ar,Kr,Xe*}, {*N,P*}, {*S,Se*}, {*Cl,Br*}, {*O,F*}, {*B,Si,C*}, {*Ti,V*}, {*As,Sb,Te*}, {*Zn,Cd,In*}, {*Hg,Tl,Pb,Bi*}, {*Cr,Mn,Fe,Co,Ni*}, {*Zr,Hf*}, {*Nb,Ta,W,Mo,Re*}, {*Cu,Ag,Au*}, {*Tc,Ru,Os,Ir*}, {*Rh,Pd,Pt*}, {*Li,Na,K,Rb,Cs*}, {*Be,Al*}, {*Mg,Ca,Sr,Ba*}, and {*Sc,Y,La*}. These results show, besides some well-known groups of similar elements, evidence of new similarity patterns, such as the diagonal relationships. A description of these patterns is the following:

*Singularity Principle* (Rodgers, 1995; Rayner-Canham and Overton, 2002): the chemistry of the second period elements often differs from that of the later members of their respective groups.

*Diagonal Relationships* (Rodgers, 1995; Rayner-Canham and Overton, 2002): there are similarities in chemical properties between an element and that at the lower right of it in the Periodic Table.

*Inert Pair Effect* (Rodgers, 1995; Rayner-Canham and Overton, 2002): in some groups, the elements following the fifth and sixth periods exhibit oxidation states two values below the maximum of their respective groups.

*Knight's move* (Rayner-Canham, 2000; Rayner-Canham and Overton, 2002; Laing, 2001): there are similarities between an element of group *n* and period *m* with the element in group *n*+2 and period *m*+1 in the same oxidation state.

*Secondary periodicity* (Ostrovsky, 2001): there are similarities between the properties of the corresponding elements belonging to period *m* and those belonging to period *m*+2.

Some other trends appear explained in Rayner-Canham 2000, Rayner-Canham and Overton 2002 and Greenwood and Earnshaw 2002.

However, the most important facts of these attempts (Zhou et al., 2000; Sneath, 2000) are 1) the use of mathematical tools to study the chemical elements and 2) the use of the properties to define each chemical element. In general, the advantage of using cluster analysis is the definition of similarity in mathematical terms but it raises the problem of selecting the number and type of properties that can be used in the clustering study.

It is important to remark that the clustering studies mentioned above do not take part of the quantum chemical studies of the elements; they are purely phenomenological ones due to the fact that they use the experimental information of the elements. These researches may be seen as a modern view of Mendeleev's methodology, taking advantage of the experimental properties to study the similarities among the elements. Finally, Restrepo et al., following this way of studying the chemical elements and the Periodic Law, developed a mathematical study of 72 chemical elements (Z=1-86, omitting Z=58-71) initially using 31 physico-chemical and chemical properties (Restrepo et al., 2004a, 2005a, 2006a). The results showed several well-known groups of the periodic table, evidencing patterns, such as: the singularity principle, the diagonal relationships, the inert pair effect and the knight's move. But an important fact of this study was the use of set theory and general topology[23], to extract and analyse[24] the results of the initial cluster analysis[25]. This mathematical approach nicely identifies the alkali metals and the noble gases as groups of elements having any similarity with other elements. Further Restrepo et al. found that the set of metals and non-metals have a particular mathematical boundary. This boundary is made of the semimetals of the lower part of the conventional

Periodic Table (Restrepo et al., 2004a, 2006a). However, in a recent study of the same set of chemical elements (Restrepo et al., 2006b) now, using 128 properties, with 103 of them being chemical, again the same results were found for the alkali metals and noble gases. On the other hand, the boundary of the set of metals was found to be formed by the elements that are considered semimetals by chemists (Rayner-Canham and Overton, 2002), from *B* to *Po* with a well-known stair-shape dividing metals from non-metals. Perhaps, the most relevant issue of these last researches is the use of mathematical tools to study the similarity among the chemical elements; a methodology similar to the one employed by Mendeleev (Mendeleev, 1869). Particularly, the last two studies show the mathematical aspect of the Periodic Law manifesting a topological mathematical structure (Villaveces, 2000). In spite of the results shown by the combination of cluster analysis and topology, this procedure has a shortcoming typical of cluster analysis – namely, the (ambiguity of) selection of the properties to make the analysis. These results (Restrepo et al., 2004a, 2004b, 2006a, 2006b) show that the selection of 31 properties (Restrepo et al., 2004a, 2006a) produces similar results to the ones obtained using 128 properties (Restrepo et al., 2006b). This raises the question regarding the number and type of needed and sufficient properties to properly define the chemical elements.

We mentioned above that nowadays we have more than 700 periodic tables and that many of them may be viewed as *meta*-descriptions of the first Periodic Table developed by Mendeleev. We can say, regarding this that all the works discussed here with the intention of finding similarities among the chemical elements, starting either from the atomic (quantum-theoretical approaches) or from the phenomenological point of view are in some sense not *meta*-descriptions of Mendeleev's Periodic Table since they are not based on the original Mendeleev Periodic Table and they do not have as their goal a better way to show the trends in the elements. The aims of these studies are the relationships in themselves. We show in figure 4 a graphical representation of the *meta*- and not *meta*-descriptions of the chemical elements.

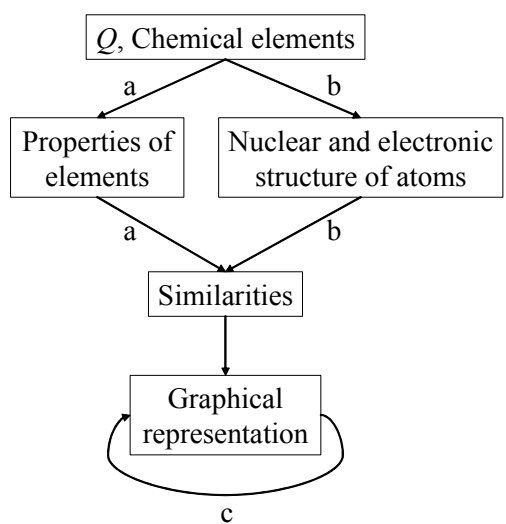

**Figure 4.** First in a) and b) the organisational philosophy for the not *meta*-descriptions and second in c) for the *meta*-descriptions; a) and b) are the phenomenological and the quantum-theoretic descriptions, respectively.

## 5. ON THE CARDINALITY OF PERIODS

Any of the mathematical studies mentioned above address mathematical explanations for the length of periods (cardinalities). The most recent attempt trying to explain this was made by Restrepo et al. (Restrepo et al., 2005b) taking into account the chemotopological approach and the number and population of the clusters found[26]. But it is possible to do a pure mathematical study regarding the cardinality of periods if we consider each period as a member of a sequence of integers.

The cardinalities of the first periods are 2, 8, 8, 18, 18, 32, 32,... An expression, which is related to these numbers is the well-known $2n^2$, where $n$ is an *integer* number. This formula does not indicate the cardinality of the periods but the maximum number of electrons quantum-theoretically allowed in each shell or electronic level $n$ (Table 1)[27].

| $n$ | 1 | 2 | 3 | 4 | ... |
|-----|---|---|----|----|-----|
| $2n^2$ | 2 | 8 | 18 | 32 | ... |

**Table 1.** Number of electrons allowed in each electronic level *n*.

Numbers produced by $2n^2$ do appear in the sequence of cardinals of the Periodic Law 2,8,8,18,18,... but they do not appear with a matching redundancy[28]. It might be interesting to develop a mathematical expression that shows these cardinals in the exact order with the proper redundancy. Regarding mathematical expressions of this sort, Sloane has collected several formulas that reproduce integer sequences and has a site on the Internet to give a compendium of these mathematical expressions (Sloane, 2006). But, our Mendeleevian sequence 2,8,8,18,18,... did not appear there[29]. Considering this, we developed a mathematical expression for these cardinals and sent it to Sloane's web site, to receive the ID number A093907. In the following we show our mathematical expression:

Let $c_n$ be the cardinality of the period $n$ in the conventional Periodic Table (Fernelius and Powell, 1982). Thus, $c_1$ is the cardinality of period 1, $c_2$ that of period 2 and so on. Then, $c_n \in Z$, where $Z$ is the set of positive integer numbers. Thus, we have $c_n = 2 \lfloor (n+2)/2 \rfloor^2$, $n \in Z$, where $\lfloor (n+2)/2 \rfloor$ represents the integer part of $(n+2)/2$.

Now we have a mathematical expression to represent the cardinality of periods of the Periodic Law[30]. This formula starts with *n=1*, but which is the end? What is the last period? The sequence shown by this formula does not end. However, the number of elements evidently is limited due to instabilities in the nuclei of atoms, as have been considered in some detail by numerous authors, and has been summarised by Karol (Karol, 2006). Indeed, there is believed to be an island of relative stability in the neighbourhood of atomic numbers between 120 and 130. Thus, according to our inductive expression for $c_n$, we would have 8 periods (Table 2).

| $n$ | 1 | 2 | 3 | 4 | 5 | 6 | 7 | 8 |
|-----|---|---|---|---|---|---|---|---|
| $2 \lfloor (n+2)/2 \rfloor^2$ | 2 | 8 | 8 | 18 | 18 | 32 | 32 | 50 |
| Accumulated number of elements | 2 | 10 | 18 | 36 | 54 | 86 | 118 | 168 |

**Table 2.** Cardinality of periods in the Periodic Law.

If we study the numbers of our sequence, they show that every $c_n$ is an even number. Taking this into account, we can define a new sequence 1,4,4,9,9,16,16,25,25,..., (half of the values in $c_n$) and we can make a graphic representation of its elements (Figure 5).

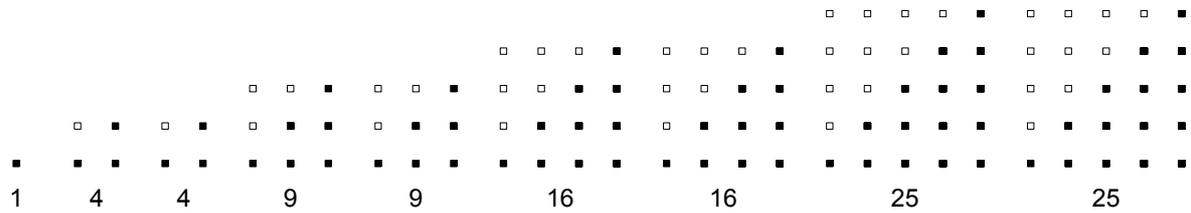

1     4     4     9     9    16    16    25    25

**Figure 5.** Graphic representation of the sequence of square numbers present in the cardinality of the periods of the Periodic Law.

Clearly this sequence is comprised from square numbers {1,4,9,16,25,...}, each repeated twice except the first. Besides, square numbers may be obtained as sums of adjacent couples of the Pythagorean triangular ($\tau\rho\iota\gamma\omega\nu o\varsigma$) numbers (Gow, 1923) (Table 3). Graphically, this may be seen in figure 5, where black points represent the triangular numbers needed to obtain the corresponding square number (black and white points)[31].

| Square numbers | 1 | 4 | 9 | 16 | 25 | ... |
|---|---|---|---|---|---|---|
| Triangular numbers | 1 | 3 | 6 | 10 | 15 | ... |

**Table 3.** Relationships between triangular and square numbers.

The sequence of triangular numbers can be defined as (Rouse, 1919; Gow, 1923): $t_n = n(n+1)/2$, $n \in Z$. Thus, we can represent the periodicity in a 3-dimensional way as is shown in figure 6a, or in a planar way as it appears in figure 6b. Figure 6a is similar to the one shown by Sugathan and Menon (Sugathan and Menon, 1956; Mazurs, 1974).

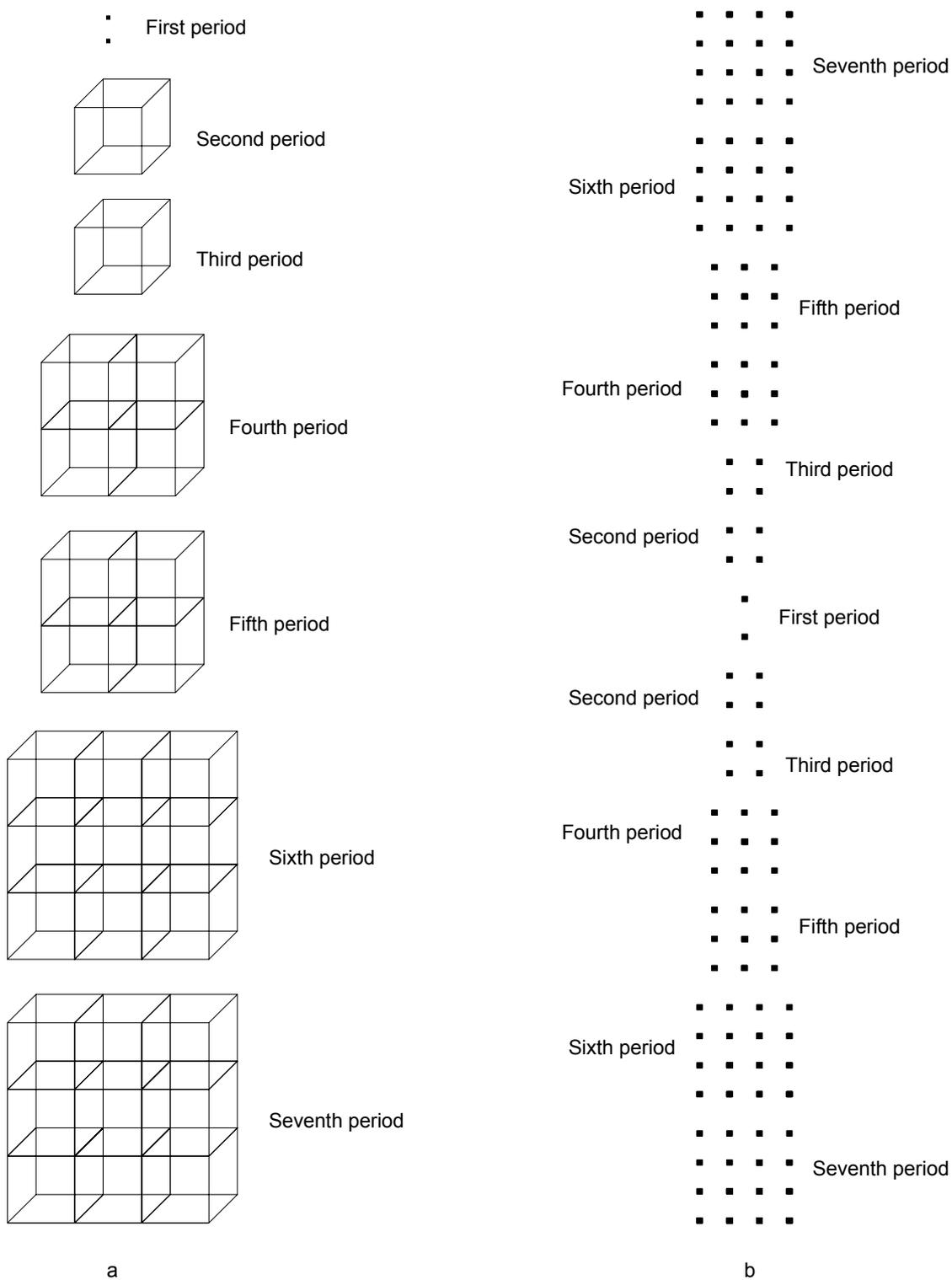

**Figure 6.** Representations of the Periodic Law based on square numbers.

Thus, it is seen that apart from the different mathematical structures of the chemical elements discussed in this paper, there is a numerical relationship among the cardinalities of the periods in the Periodic Law. The triangular numbers discovered by Pythagoras can express these numerical relations[32].

## 6. CONCLUSIONS

We have indicated different approaches developed to study the Periodic Law and the chemical elements from the mathematical point of view. We mentioned the first approaches developed during the end of the 19[th] century until those done in the 21[st] century. The different mathematical theories and tools used to study the Periodic Law can be classified in the following groups: numerical attempts, information theoretical approaches, order theoretical studies, similarity and topological approaches. Several of the researches mentioned provide mathematical structures for the chemical elements. These structures are: order-theoretic, topological, and numerical structures. It is important to note that despite the quantum chemical studies, very common nowadays and briefly mentioned in this paper, there are some others related to the phenomenological character of the chemical elements. One important conclusion regarding all the researches showed in this paper is the evidence of the current importance of the Periodic Law in our times. This is notorious since novel mathematical theories have been applied to the study of the chemical elements as can be seen in the application of the information theory, set theory, partial ordered sets, and topology. In addition, it is possible to find connections with some ancient philosophical teaching, such as the Pythagorean doctrine. Some authors have thought that the beauty in Nature is a fact and a way to express this fact is through mathematics. This paper evidences such 'beauty' of the chemical elements.

According to ancient Greek thinking Nature has a principle of homogeneity. This fundamental homogeneity was studied from two points of view: the unity in the matter, and the unity in the structure of things. Particularly, Pythagoras pursued the latter way (Gow, 1923) based on what he learned about numbers in Egypt. He realised that numbers were fundamental to the exact description of forms and their relationships, and concluded that number was the cause of form and every other quality. In this way he pointed out 'Number is quantity and quantity is form and form is quality.' (Gow, 1923). This is the foundation of what we know regarding the Pythagorean doctrine that considers number as the basis of creation (Gow, 1923; Russell, 1993). Taking advantage of this and assuming that the chemical elements are the forerunner of the whole of the materials as a premise, we can then say that this philosophical doctrine has not lost its use even after twenty-five centuries as the Periodic Law shows an underlying Pythagorean structure.

## NOTES

1. Mendeleev predicted the following properties, particularly to *Ge*: atomic weight, density, molar volume, melting point, specific heat, valence, colour, the method of isolation, reactivity with *HCl* and *NaOH*, reactivity with oxygen, empirical formula of the compound with oxygen ($GeO_2$), density of the compound $GeO_2$, solubility of $GeO_2$ in acid, solubility of $GeS_2$ in water and ammonium sulfide, boiling point of $GeCl_4$ and boiling point of $GeEt_4$ (Greenwood and Earnshaw, 2002). Mendeleev however called the element '*eka*-silicon' as his predictions were in advance of experimental evidence of the element. Indeed, Mendeleev's work gained prominence when it turned out that he made significantly more accurate predictions of the density of $GeO_2$ than the initially (later) reported experimental density which turned out to be in error.

2. Mathematical structure here indicates a set with a relation among its objects. Formally we can say: an ordered pair ($Q,r$) is called a structure if $r$ is a relation on $Q$ (Potter, 2004).

3. Cardinality coincides for finite sets with the number of elements of a set.

4. It is important to note that this sort of periodicity does not correspond to the mathematical meaning of singly periodic. A property $f(x)$ is singly periodic if $f(x) = f(x + np)$, where $n$ is the period. Examples of singly periodic functions, where $n = 2\pi$, are $f(x) = sin(x)$ and $g(x) = cos(x)$, for instance.

5. However, Restrepo et al. have shown (Restrepo et al., 2004a, 2006a) that the similarity relationships of the chemical elements are not dependent on the atomic number or the atomic weight.

6. We consider Lewis' work as mathematical. Madelung's work and all the following researches done from quantum chemistry also may be considered as mathematical considering that they arise from the solution of an *eigenvalue* equation $H\Psi=E\Psi$. But our interest in this paper is to consider the researches different from the ones that arise from the quantum mechanics.

7. There are two general opinions in the literature regarding the authorship of the Madelung Rule. One of them indicates that this rule was proposed by Bohr in 1922 (Scerri et al., 1998), though according to Ostrovsky (Ostrovsky, 2001) Bohr's paper cited by Scerri et al. (Scerri et al., 1998) does not have a formulation of this rule. The other general opinion credits the rule to Madelung in his handbook of 1936 (Madelung, 1936). According to Ostrovsky (Ostrovsky, 2001): "It is difficult to trace the origin of this rule that looks a kind of scientific folklore." Nowadays a wide-spread idea is that Bohr was the one who proposed the rule e.g. Löwdin (Löwdin, 1969), who ascribes the rule to Bohr without giving references.

8. A comprehensive overview of the researches regarding the Periodic Law from a quantum chemical view has been made by Scerri et al. (Scerri et al., 1998).

9. According to Bonchev, "some predictions were made for a number of properties of the transactinide chemical elements 113-120 and for the nuclear binding energies of the nuclides of the elements 101 to 108" (Bonchev, 2006).

10. It is possible to define a partial order on a set $A$ as a binary relation $>$ satisfying the following conditions (Lipschutz, 1965):

   1) $a \in A \Rightarrow a > a$

   2) $a > b \wedge b > a, a,b \in A \Rightarrow a = b$

   3) $a > b \wedge b > c, a,b,c \in A \Rightarrow a > c$

   Here $>$ may for instance mean 'exceeds', 'greater than' or 'older than'.

11. Regarding the order in the chemical elements, another important aspect of Mendeleev's work was not only the study of the similarity relationships of the elements but the application of the notion of order found in the Periodic Law. It was possible to perceive the Law when Mendeleev organised the atomic weights of the elements and plotted his graphic against the valence. Mendeleev's formulation entailed a correspondence between the set of atomic weights and the set of natural numbers. That is, it entailed the mathematical concept of order.

12. Total order means that for every pair of elements $a, b \in A$, it is always possible to have $a > b$ or $b > a$. An example of total order is the set of real numbers with the natural order defined by $a \leq b$.

13. Here it is possible to establish a connection between order and topology. Two shells or, in general, two objects $a, b$ with the relation $>$ can be considered as neighbours, where one neighbour is greater than the other, for instance. In this way we can talk about binary directed neighbourhoods governed by the order relation $>$.

14. For every $n, n', l, l', a, b \in \mathrm{N}$, $(n, l) < (n', l')$ iff $(n+a, l+b) < (n'+a, l'+b)$.

15. However, the Periodic Table of ions is different than that of the neutral elements as has been suggested by Pettifor (Pettifor, 1984, 1986, 1988a, 1988b, 1996) and recently shown by Railsback (Railsback, 2003). In terms of Klein's ideas (Klein, 1995, 2000; Klein and Babić, 1997), the Periodic Table of ions should manifest a multi-poset structure different from the one found in the conventional Periodic Table.

16. It is important to remark that this order is different from the one of Kreinovich et al. since the Pettifor order is phenomenological, being based on experimental properties. Kreinovich et al. order is based on the shells of the atoms.

17. Pettifor order is: He,Ne,Ar,Kr,Xe,Rn,Fr,Cs,Rb,K,Na,Li,Ra,Ba,Sr,Ca,Yb,Eu,Sc,Lu,Tm,Er,Ho, Dy,...,Tb,Gd,Sm,Pm,Nd,Pr,Ce,La,Lr,No,Md,Fm,Es,Cf,Bk,Cm,Am,Pu,Np,U,Pa,Th,Ac,Zr,Hf ,Ti,Ta,Nb,V,W,Mo,Cr,Re,Tc,Mn,Fe,Ru,Os,Co,Rh,Ir,Ni,Pt,Pd,Au,Ag,Cu,Mg,Hg,Cd,Zn,Be,T l,In,Al,Ga,Pb,Sn,Ge,Si,B,Bi,Sb,As,P,Po,Te,Se,S,C,At,I,Br,Cl,N,O,F,H.

18. Each property was viewed as a function $f(w_i)$ of the atomic weight $w_i$. According to Mendeleev (Mendeleev, 1869) the fact of having similar properties for some elements "at once raises the question whether the properties of the elements are expressed by their atomic weights and whether a system can be based on them."

19. Particularly in the case of chemical elements, if we define each element using 7 of its properties, then we have a 7-dimensional space, and each element is mathematically represented as an ordered 7-tuple.

20. Fuzzy cluster analysis is a special case of cluster analysis where there are no binary relationships of belonging. Then, it is not necessary that an element *a* either belongs to or does not belong to a set. Fuzzy cluster analysis establishes that it is possible to consider that the element *a* belongs to 'in part' to a set.

21. Sneath and Sokal are the founders of the Numerical Taxonomy (Sokal and Sneath, 1963), nowadays called cluster analysis. This tool was first used by them to classify biological species taking into account morphological properties without any consideration of phylogenetic relationships. Here, 'similarity' arises from comparison of the properties not from a priori classifications (Felsenstein, 2004).

22. Sneath did not consider more elements due to the lack of experimental data to define the rest of the elements. The problem of the necessary number of properties to define the chemical elements is still an open question. Regarding this, Mendeleev said (Mendeleev, 1869): "The numerical data for simple bodies are limited at the present time." It would be very interesting to have a defined set of properties to define the chemical elements as it is possible to have a set of three properties to define the electron nowadays, {mass, charge, spin} (Villaveces, 2004).

23. This new methodology that combines chemometrics and topology was called chemotopology (Restrepo et al., 2005a).

24. Topology is a mathematical theory whose history dates back at least to 1736 when Euler presented a solution (Euler, 1741, date of publication) to the problem of finding a path over the seven bridges of Königsberg (Germany) in such a way that each bridge was crossed just once (Biggs et al., 1998). Particularly, these sorts of problems are now studied by the graph theory, one branch of topology (Harju, 2002). But it was Listing (Listing, 1847), in the 19[th] century, who used the word topology for the first time in a book titled *Vorstudien zur Topologie* (Preliminary studies on topology). Such a word condenses the seminal idea of Euler about the study of the place ($τόπος$: place, $λόγος$: thought). Euler in his paper on the Königsberg bridges (Euler, 1741), titled *Solvtio problematis ad geometriam sitvs pertinentis* (The solution of a problem relating to the geometry of position), says the following (the original Latin paper appears in reference Euler, 1741):

In addition to that branch of geometry which is concerned with magnitudes, and which has always received the greatest attention, there is another branch, previously almost unknown, which *Leibniz* first mentioned, calling it *geometry of position*. This branch is concerned only with the determination of position and its properties; it does not involve measurements, nor calculations made with them. It has not yet been satisfactory determined what kinds of problems are relevant to this geometry of position, or what methods should be used in solving them. Hence, when a problem was recently mentioned, which seemed geometrical but was so constructed that it did not require the measurement of distances, nor did calculation help at all, I had no doubt that it was concerned with the geometry of position – specially as its solution involved only position, and no calculation was of any use. I have therefore decided to give here the method which I have found for solving this kind of problem, as an example of the geometry of position (Biggs et al., 1998).

25. We can use topology when we have a set and the relationships among the elements belonging to the set. In the first example of Euler (Euler, 1741) he had a set of seven bridges and he knew their relationships. In the case of chemical elements we have something similar to Euler's problem. We have a set $Q$ of chemical elements and the similarity relationships among them. The similarity relationships, extracted from the results of cluster analysis, are condensed in $\tau$, and the couple $(Q, \tau)$ is called a (discrete) topological space. Thus, the chemical elements and their similarities are a topological space. The mathematical development of this methodology (chemotopology) appears in the references Restrepo et al., 2004a, 2006a and 2006b.

26. This method is based on the number $N$ of clusters that appear when one cuts branches on the tree (dendrogram) and when one considers the population of each branch. The method establishes that the clusters that offer more neighbourhood information are those that arise from a cut that gives $k<N$ clusters and a value $\prod_{i=1}^{N} p_i$ of the population of the clusters. Here $p_i$ is the population of the cluster $i$. A description of this methodology appears in references (Restrepo et al., 2005b, 2005c; Uribe et al., 2005).

27. Regarding this, Pauling said: "The successive electron shells [of the noble gases ...] involve the following numbers of electrons: 2, 8, 8, 158, 18, 32, 32. These numbers are equal to the numbers of elements in the successive periods of the periodic system." (Pauling, 1970).

28. Here $2n^2$ gives the set of numbers that appears in the sequence of cardinals, but not the sequence. Thus, if we called $C=\{x \mid x=2n^2,\ n$ is a positive integer$\}$ then we have that $C=\{2,8,18,32,...\}$, but not the ordered set $(2,8,8,18,18,32,32,...)$.

29. However, this site does include several sequences of chemical phenomena such as the sequence of the numbers of alkanes, initially worked by Cayley in 1875. A compendium of Cayley's papers appear in Cayley, 1896.

30. Weise (Weise, 2003) expresses the total number of electrons of noble gases atoms by the expression $Z_n=((-1)^n(3n+6)+2n^3+12n^2+25n-6)/12$, where $n$ is an *integer number* ($n=1,2,\ldots$). He also developed a method for getting the total number of electrons of the noble gases by modifications of the Pascal's triangle. Despite the success for reproducing the sequence of total electrons in noble gases, Weise did not develop a mathematical expression for the sequence of the cardinalities of the periods in the Periodic Law.

31. Weise (Weise, 2003), through modifications of Pascal's triangle, used the triangular numbers for reproducing the sequence of total electrons of the noble gases.

32. Pythagoras classified all numbers as 'odd' ($\alpha\rho\tau\iota o\iota$) or 'even' ($\pi\epsilon\rho\iota\sigma\sigma o\iota$). The odd numbers were also called 'gnomons' ($\gamma\nu\omega\mu o\nu\epsilon\varsigma$) and their sum from 1 to $(2n+1)$ was called a 'square' ($\tau\epsilon\tau\rho\alpha\gamma\omega\nu o\varsigma$) (Gow, 1923).

**ACKNOWLEDGEMENTS**

aspects of the mathematical structures of chemical sets. A thank is also made to Dr. R. Hefferlin from the Southern Adventist University (USA) for having provide information and also the manuscript by D. Weise. Finally, the Universidad de Pamplona (Colombia) and specially Dr. A. González, rector of that University, are thanked for their support to develop research on mathematical and philosophical chemistry.

# REFERENCES

M. S. Aldendefer and R. K. Blashfield. *Cluster Analysis (Sage University Paper series on Quantitative Applications in the Social Sciences, N°. 44)*. Sage University: Newbury Park, 1984.

E. V. Babaev and R. Hefferlin. *The Concepts of Periodicity and Hyperperiodicity: from Atoms to Molecules*. In Concepts in Chemistry; D. H. Rouvray and E. R. Kirby Eds.; Research Studies Press: Taunton, 1996.

N. L. Biggs, E. K. Lloyd and R. J. Wilson. *Graph Theory 1736-1936*. Clarendon Press: Oxford, 1998.

N. Bohr, Der Bau der Atome und die physikalishen und chemishen Eigenschaften der Elemente. *Zeitschrift für Physik*, 9: 1 67, 1922. English translation in *N. Bohr, Collected Works*, J. R. Nielsen, Ed.; North-Holland: Amsterdam, 1977.

D. Bonchev. *Periodicity of the Chemical Elements and Nuclides: An Information Theoretic Analysis*. In The Mathematics of the Periodic Table; R. B. King and D. H. Rouvray Eds.; Nova: New York, 2006.

R. Carbó, L. Leyda and M. Arnau. How Similar is a Molecule to Another? An Electron Density Measure of Similarity between Two Molecular Structures. *Int. J. Quantum Chem*. 17, 1185-1189, 1980.

R. Carbó and L. Domingo. LCAO-MO Similarity Measures and Taxonomy. *Int. J. Quantum Chem*. 22, 517-545, 1987.

A. Cayley. *The collected mathematical papers of Arthur Cayley*. Vol IX, 1st ed. Cambridge University Press: Cambridge, 1896. Internet: http://www.hti.umich.edu/cgi/t/text/pageviewer-idx?c=umhistmath;cc=umhistmath;sid=0128b870e72ce32c3182761dac22f43e;q1=triangular%20numbers;rgn=full%20text;view=pdf;seq=00000001;idno=ABS3153.0009.001

D. Cruz-Farritz, J. A. Chamizo and A. Garritz. *Estructura atómica, un enfoque químico*. Addison-Wesley: Wilmington, 1991.

L. Euler. Solvtio Problematis ad Geometriam Sitvs Pertinentis. *Commentarii Academiae Scientiarum Imperialis Petropolitanae* 8, 128-140, 1741. *Opera Omnia*. Series 1. Volume 7, pp. 1-10. Internet: http://www.eulerarchive.com/

B. S. Everitt. *Cluster Analysis*. Arnold: London, 1995.

J. Felsenstein. *Inferring Phylogenies*. Sinauer: Sunderland, 2004.

W. C. Fernelius and W. H. Powell. Confusion in the Periodic Table of the Elements. *J. Chem. Educ.* 59, 6, 504-508, 1982.

F. L. G. Frege. *Grundgesetze der Arithmetik*. Pohle, Band I: Jena, 1893.

F. L. G. Frege. *Grundgesetze der Arithmetik*. Pohle, Band II: Jena, 1903.

A. D. Gordon. *Classification*. Chapman and Hall: London, 1981.

G. Gow. *A short History of Greek mathematics*. C. E. Stechert: New York, 1923.

N. N. Greenwood and A. Earnshaw. *Chemistry of the Elements.* Butterworth Heinemann: Oxford, 2002.

T. Harju. *Lecture Notes on Graph Theory*. Department of Mathematics University of Turku: Turku, 2002.

P. J. Karol. *The Heavy Elements*. In The Mathematics of the Periodic Table; R. B. King and D. H. Rouvray Eds.; Nova: New York, 2006.

D. J. Klein. Similarity and dissimilarity in posets. *J. Math. Chem*. 18, 321-348, 1995.

D. J. Klein and D. Babić. Partial Orderings in Chemistry. *J. Chem. Inf. Comput. Sci.* 37, 656-671, 1997.

D. J. Klein. Prolegomenon on Partial Orderings in Chemistry. *MATCH*. 42, 7-21, 2000.

M. Laing. Periodic Patterns. *J. Chem. Educ*. 78, 7, 877, 2001.

G. N. Lewis, The Atom and the Molecule. *J. Am. Chem. Soc*., 38: 762-786, 1916. Internet: http://dbhs.wvusd.k12.ca.us/webdocs/Chem-History/Lewis-1916/Lewis-1916.html

S. Lipschutz. *General Topology*. McGraw-Hill: New York, 1965.

J. B. Listing. *Vorstudien zur Topologie*. Gottinger Studien, 1847.

P. O. Löwdin. *International Journal of Quantum Chemistry (Symposium)* IIIS: 331-334, 1969.

E. Madelung, *Die mathematischen Hilfsmittel des Physiker*. Springer: Berlin, 1936.

E. G. Mazurs. *Graphic Representations of the Periodic System During One Hundred Years*. University Alabama Press: Alabama, 1974.

D. I. Mendeleev, Sootnoshenie svoistv s atomnym vesom elementov. *Zh. Russ. Khim. Obshch* 1 (2/3): 60-77, 1869.

D. I. Mendeleev, The Periodic Law of the Chemical Elements. *J. Chem. Soc.,* 55: 634-656, 1889. Internet: http://web.lemoyne.edu/~giunta/mendel.html

L. Meyer, Die Natur der chemischen Elemente als Function ihrer Atomgewichte. *Annalen der Chemie und Pharmarcie Supplementband* 7: 354–364, 1870.

H. G. J. Moseley, The High Frequency Spectra of the Elements. *Phil. Mag*, 1024, 1913. Internet: http://dbhs.wvusd.k12.ca.us/webdocs/Chem-History/Moseley-article.html

V. N. Ostrovsky. What and How Physics Contributes to Understanding the Periodic Law. *Found. Chem*. 3, 145-182, 2001.

L. Pauling. *General Chemistry*. Dover: New York, 1970.

D. G. Pettifor. A chemical scale for crystal-structure maps. *Sol. State Commun*. 51, 1, 31-34, 1984.

D. G. Pettifor. The structures of binary compounds: I. Phenomenological structure maps. *J. Phys. C: Solid State Phys*. 19, 285-313, 1986.

D. G. Pettifor. Structure maps for pseudobinary and ternary phases. *Mat. Sci. Tech*. 4, 675-692, 1988a.

D. G. Pettifor. Structure maps in magnetic alloy design. *Physica B*. 149, 3-10, 1988b.

D. G. Pettifor. Phenomenology and Theory in Structural Prediction. *J. Phase Equil*. 17, 5, 384-395, 1996.

M. Potter. *Set Theory and its Philosophy*. Oxford: Oxford, 2004.

L. B. Railsback. An earth scientist's periodic table of the elements and their ions. *Geology*. 31, 9, 737-740, 2003.

G. Rayner-Canham. Periodic Patterns. *J. Chem. Educ*. 77, 8, 1053-1056, 2000.

G. Rayner-Canham and T. Overton. *Descriptive Inorganic Chemistry*. Freeman: New York, 2002.

G. Restrepo, H. Mesa, E. Llanos and J. L. Villaveces. Topological Study of the Periodic System. *J. Chem. Inf. Comput. Sci*. 44, 68-75, 2004a.

G. Restrepo, E. J. Llanos and H. Mesa. 2004, *Chemical Elements: A Topological Approach*, Theodore Simos, *Proceedings of the International Conference of Computational Methods in Sciences and Engineering 2004*. 753 – 755, Athens, Greece, 19 – 23 November 2004. Utrecht: VSP, 2004b.

G. Restrepo and J. L. Villaveces. From trees (dendrograms and consensus trees) to topology. *Croat. Chem. Acta* 78, 275-281, 2005a.

G. Restrepo, E. J. Llanos and J. L. Villaveces. 2005, *Trees (dendrograms and consensus trees) and their topological information*, Subhash C. Basak and Indira Ghosh, *Proceedings of the Fourth Indo-US Workshop on Matehmatical Chemistry*. 39 – 62, Pune, India, 8 – 12 January 2005. Pune, 2005b.

G. Restrepo and R. Brüggemann. Ranking regions through cluster analysis and posets. *WSEAS Trans. Inf. Sci. Appl*. 2, 976-981, 2005c.

G. Restrepo, H. Mesa, E. J. Llanos and J. L. Villaveces. *Topological Study of the Periodic System*. In The Mathematics of the Periodic Table; R. B. King and D. H. Rouvray Eds.; Nova: New York, 2006a.

G. Restrepo, E. J. Llanos and H. Mesa. Topological space of the chemical elements and its properties. *J. Math. Chem.* 39, 401-416, 2006b.

G. Restrepo, H. Mesa and J. L. Villaveces. On the Topological Sense of Chemical Sets. *J. Math. Chem.* 39, 363-376, 2006c.

D. Robert and R. Carbó-Dorca. On the extension of quantum similarity to atomic nuclei: Nuclear quantum similarity. *J. Math. Chem.* 23, 327-351, 1998.

D. Robert and R. Carbó-Dorca. General Trends in Atomic and Nuclear Quantum Similarity Measures. *Int. J. Quantum Chem.* 77, 685-692, 2000.

G. E. Rodgers. *Química Inorgánica*. McGraw-Hill: Madrid, 1995.

W. W. Rouse. *A short account of the history of mathematics*. Macmillan and co: London, 1919. Internet: http://www.hti.umich.edu/cgi/t/text/pageviewer-idx?c=umhistmath;cc=umhistmath;sid=0128b870e72ce32c3182761dac22f43e;q1=triangular%20numbers;rgn=full%20text;view=pdf;seq=00000001;idno=ACA1117.0001.001

B. Russell. *Introduction to Mathematical Philosophy*. Dover: New York, 1993.

E. Scerri, V. Kreinovich, P. Wojciechowski and R. Yager. Ordinal explanation of the periodic system of chemical elements. *Int. J. Uncertainty, Fuzziness and Knowledge-Based Systems*, 6, 387-400, 1998.

N. Sloane. *The On-Line Encyclopedia of Integer Sequences*. Internet: http://www.research.att.com/~njas/sequences/

P. H. A. Sneath. Numerical Classification of the Chemical Elements and its relation to the Periodic System. *Found. Chem.* 2, 237-263, 2000.

R. R. Sokal and P. H. A. Sneath. *Principles of Numerical Taxonomy*. W. H. Freeman: San Francisco, 1963.

K. K. Sugathan and T. C. K. Menon. Periodic classification and electronic configuration of elements. *Current Sci.* 25, 85, 1956.

A. P. Sutton. *Electronic Structure of Materials*. Oxford Science Publications: Oxford, 1996.

E. A. Uribe, M. C. Daza and G. Restrepo. Chemotopological study of the fourth period monohydrides. *WSEAS Trans. Inf. Sci. Appl.* 2, 1085-1090, 2005.

J. L. Villaveces. Química y Epistemología: una relación esquiva. *Revista Colombiana de Filosofía de la Ciencia* 1, 2-3, 9-26, 2000.

J. L. Villaveces. Personal communication, 2004.

D. Weise. *A Pythagorean approach to problem of periodicity in chemical and nuclear physics*. In Advanced topics in theoretical chemical physics; J. Maruani, R. Lefebvre and E. J. Brändas Eds.; Kluwer: Dordrecht, 2003.

X. Z. Zhou, K. H. Wei, G. Q. Chen, Z. X. Fan and J. J. Zhan. Fuzzy Cluster Analysis of Chemical Elements. *Jisuanji Yu Yingyong Huaxue*. 17, 167-168, 2000.